\documentclass[11pt, fleqn]{article}
\usepackage{amsmath}
\usepackage{amssymb,amsthm,epsfig}
\usepackage{enumerate}
\usepackage{mathabx}
\usepackage{statmath}
\usepackage{accents}
 \def\vt{t\kern-0.22em\raise.18ex\hbox{\char'47}\lower.18ex\hbox{}\kern-0.08em}
\newtheorem{theorem}{Theorem}[section]

\newtheorem{ob}{Observation}[section]

\newtheorem{rem}{Remark}[section]

\newcommand{\old}[1]{{}} 

\usepackage{accents}
\newlength{\dhatheight}

\newcounter{obr}
\newcounter{tabul}
\begin{document}
\title{Bounding the edge cover  of a hypergraph  \
}
\author{Farhad Shahrokhi\\
Department of Computer Science and Engineering,  UNT\\
P.O.Box 13886, Denton, TX 76203-3886, USA\
Farhad.Shahrokhi@unt.edu
}

\date{}
\maketitle
\date{} \maketitle


\begin{abstract}
 Let $H=(V,E)$ be a hypergraph.  Let $C\subseteq E$, then $C$ is an {\it edge cover}, or a {\it set cover}, if $\cup_{e\in C} \{v|v\in e\}=V$. A subset of vertices $X$ is {\it independent}  in $H,$  if no two vertices in $X$ are in any edge.  Let  $c(H)$ and   $\alpha(H)$  denote the cardinalities  of a smallest edge cover and largest  independent set in $H$, respectively. 
We show that $c(H)\le {\hat m}(H)c(H)$,   where ${\hat m}(H)$ is a parameter called the {\it mighty  
degeneracy} of $H$.  
 Furthermore, we show that the inequality is tight and demonstrate the   applications in domination theory. 
 
\end{abstract}

\section{Introduction} 
We assume the reader is familiar with standard graph theory \cite{Char}, hypergraph theory \cite{B}, \cite{Bo},  domination theory \cite{HHS}, and algorithm analysis \cite{CLR}. 
Throughout this paper we denote by $H=(V,E)$ a hypergraph on vertex set $V$ and the edge set $E$. So any $e\in E$ is a subset of $V$.  We  do not  allow multiple edges in our definition of a hypergraph, unless explicitly stated.  Every hypergraph can be represented by its incidence  bipartite graph $B$ whose vertex set is $V\cup E$. If $x\in V$ and $e\in E$, then $xe$ is an edge in $B$,  provide that $x\in e$. Let $C\subseteq E$, then $C$ is an {\it edge cover}, or a {\it set cover},   if $\cup_{e\in C} \{v|v\in e\}=V$. A subset of vertices $X$ is {\it independent}  in $H$, 
if no two vertices in $X$ are in any edge. 
 Let  $c(H)$  and $\alpha(H)$ denote the cardinalities  of a largest   independent set and a smallest edge cover  in $H$, respectively. 
 It is known that computing $\alpha(H)$ and 
$c(H)$ are NP hard problems \cite{GJ}.  
Clearly, $c(H)\ge\alpha(H)$. Furthermore, it is known that $c(H)$ can not bounded above by a function of $\alpha(H)$, only. 
However, an  important result in this area is known. Specifically,  it is a consequence of the result   in \cite{Sy} that    
\begin{equation}\label{e1}
c(H)={\alpha(H)}^{O(2^v)}   
\end{equation}
where $v$ denotes the vc dimension of $H$ \cite{vc}. 
 Design of approximation algorithms for the edge cover problem has been an active and ongoing research in computer science. A greedy algorithm  \cite{ch}, \cite{lo}  is known to approximate $c(H)$ within $O(log(n)$ from its  optimal value. Moreover, there are examples of hypergraphs that show the  worst case approximation scenario of $O(log(n))$ can not be improved \cite{BG}. 

The main result of this paper is to show that  \begin{equation}\label{eq0}
c(H)\le {\hat m}(H)\alpha(H)
\end{equation} 
where the multiplicative factor  ${\hat m}(H)$ is a parameter called the {\it mighty  degeneracy } of $H$  which we introduce here.  Recall that a set $S\subseteq V$ is a  {\it transversal set}  (hitting set) in the hypergraph $H=(V,E)$, if every $e\in E$ has a vertex in $S$. A  set $M\subseteq E$ is a matching  in $H$,  if every  two edges in $M$ are disjoint. 
 Let $\tau(H)$ and $\rho(H)$ denote  the sizes of a smallest transversal and a largest matching in $H$, respectively,  and
note that $\tau(H)\ge \rho(H)$. 

 A direct consequence of (\ref{eq0}) is  that
\begin{equation}\label{eq1}
\tau(H)\le {\hat m}(H^d)\rho(H)
 \end{equation}
  where ${\hat m}(H^d)$ is the  mighty   degeneracy of the {\it dual hypergraph} of $H$, defined as $H^d=\{E,V\}$.

This paper is organized as follows. In Section Two we introduce  some terms and concepts and set up our notations. Particularly, we introduce  the   {\it strong  degeneracy}  of a hypergraph, denoted by 
 ${\hat s}(H)$,   which  is an upper bound on  ${\hat m }(H)$.   In Section Three we derive (\ref{eq0})  which is the main result, and  also present a linear time  algorithm for computing ${\hat s}(H)$.  Section Four contains the applications to domination theory of graphs. Specifically, we show  ${\hat s}(H)=1$ (and hence ${\hat m}(H)=1$),   when the underlying graph $G$ is a tree and $H$ is the so called closed or open neighborhood hypergraph of $G$. Consequently,  we  provide   new proofs (and algorithms) for two classical results in  domination theory \cite{MM}, \cite{Rall}, by  showing that in  any tree  the size of a smallest dominating (total domination)  set equals to the size of a largest 2- packing (open 2-packing). The results in Section Four are conveniently derived  utilizing concept of strong degeneracy, instead of mighty degeneracy, however  generally speaking, the  former can be much large than the latter. In Section Five we give  examples   of hypergraphs with bounded  mighty degeneracy, whose strong  degeneracy is a linear function of number of vertices. Section Six contains our suggestions for future research.

\section{Preliminaries} 
 Let $H=(V,E)$,  let $S\subseteq V$ and $e\in E$, then $e\cap S$ is the  {\it trace} of $e$ on $S$. 
The {\it restriction} of  $H$ to  $S$, denoted by $H[S]$,  is the  hypergraph  on vertex set $S$ whose edges are the set of all distinct  traces of edges in $E$ on $S$. $H[S]$ is also referred to as the {\it induced subhypergraph}  of $H$ on $S$. In general,  a  hypergraph $I$ is a subhypergraph of $H$, 
if it can be obtained by removing some vertices and some  edges from $H$\footnote {When a vertex set is removed from $H$, the edges of $H$ will also be  updated accordingly.}. $S$ is {\it shattered} in $H$, if  any $X\subseteq S$ is a trace. Thus if $S$ is shattered, then it has $2^{|S|}$ traces.  The Vapnik–Chervonenkis (VC) dimension of a hypergraph $H$, denoted by $vc(H)$,  is the cardinality of the largest subset  of $V$ which is shattered in $H$. 
Let $H=(V,E)$ and let  $x\in V$. The {\it degree}  of $x$ denoted by $d_H(x)$ is the number of edges that contain $x$. 
The  {\it strong degree} of $x$ in $H$, denoted by  $s_H(x)$,  is  the number of  distinct {\it maximal} edges  that contain $x$. 
( An edge is maximal, if it is not properly contained in another edge.) 
Let   
$\delta(H)$ and $s(H)$ denote  the smallest degree and  smallest  strong degree, respectively, of any vertex in $H$.   The {\it degeneracy} and {\it strong degeneracy} of $H$, denoted by ${\hat\delta}(H)$ and  ${\hat s}(H)$, respectively,  are the largest minimum degree and largest minimum strong degree of any induced  subhypergraph  of $H$. 
Let $R\subseteq S$.  A {\it strong subset} of $V$ in $H$ is a non empty subset of $V$ which is  obtained by  removing all vertices in $R$ from $H$,  as well as all vertices in the edges  that have nonempty intersection with $R$  and  all vertices in such edges. 

 The {\it mighty  degeneracy} of $H$, denoted by 
${\hat m}(H)$, is the largest minimum strong degree of any strong subhypergraph of $H$.   Clearly, for any $x\in V$  one has
${s}_H(x)\le {d}_H(x)$ and consequently 
\begin{equation}\label{e2} 
 {\hat m}(H)\le {\hat s}(H)\le {\hat \delta}(H)
 \end{equation}

 

\section{Our Greedy Algorithms}

Our next result is the main result of this paper.

\begin{theorem}\label{t0}
{\sl Let $H=(V,E)$ be a hypergraph, then there is an an edge cover $C$,  and an independent set $X$ in $H$ so that 
\begin{equation}\label{e4}
  |C|\le {\hat m}(H)|X|
  \end{equation}
  Consequently 
  \begin{equation}\label{e5}
  |C|\le {\hat s}(H)|X|
  \end{equation}
  
Moreover, $X$ and and $C$ can be constructed in $O(|V|+\sum_{e\in E}|e|)$ time. }
\end{theorem}
{\bf Proof.} Consider  the following algorithm.

Initially,  set $i\leftarrow 1$,  $I\leftarrow H, W\leftarrow V$ and $K\leftarrow E$.    While there are vertices in $W$ repeat the following steps: Remove the  vertex  of minimum strong   degree, denoted by $x_i$, from $W$, remove  the set of all distinct maximal edges  containing $x_i$ from $K$, then, remove  and the set of all vertices contained in theses edges from  $W$ and finally set  $i\leftarrow i+1$. 

Clearly, the algorithm terminates. Now let $t$ be the number iterations of the algorithm and at any iteration $i=1,2,...,t$,   let  $I_i$ denote the constructed hypergraph, (which is strongly induced), and let $W_i$ (which is a strong subset) and $K_i$ denote,   respectively,  the vertices and edges of  this hypergraph. Let  $X=\{x_1,x_2,...,x_t\}$ be the set of all vertices removed from $H$ when the algorithm terminates. Clearly, $X$ is an independent set in $H$. We denote  by $K_{x_i}$ the set of all distinct maximal edges  containing the vertex $x_i$ in the hypergraph $I_i$ at iteration of $i$ of the
algorithm and note that  $|K_{x_i}|\le {\hat m}(H)$, since $x_i$ is the vertex of minimum strong degree in $I_i$.
Consequently,  
\begin{equation}\label{eq20}
\sum_{i=1}^t |K_{x_i}|\le   {\hat m}(H)\times t={\hat m}(H)\times |X|
\end{equation}  

Now for $i=1,2,...,t$,  let $C_{x_i}$ be the set of all edges in $H$ obtained by extending  each edge of  $K_{x_i}$   in $I_i$ 
to an edge in $H$ and let $C=\cup_{i=1}^t C_{x_i}$. Clearly,  $C$ is 
an edge cover and furthermore $|C|=|\cup_{i=1}^t F_{x_i}|$, and therefore the first claim follows from (\ref{eq20}).

 To verify the second inequality note that ${\hat m}(H)\le {\hat s}(H)$. We omit the details of claims regrading time complexity that involves representing  $H$ as a bipartite graph. 
$\Box$

To use Theorem \ref{t0} we really need to know ${\hat m}(H)$. Alternatively,  we can use  ${\hat s}(H)$ which is an upper bound for ${\hat m}(H)$. At this time, we still do not know how to efficiency compute ${\hat m}(H)$. We finish  this section by presenting  a simple greedy algorithm for computing ${\hat s}(H)$ which is  similar to the known algorithm for computing degeneracy of $H$, or ${\hat \delta}(H)$.  The properties of the output of algorithm will be used to prove  our results in the next section.  

\begin{theorem}\label{t1}
{\sl Let   $H=(V,E)$ be a hypergraph on $n$ vertices,  then  
${\hat s}(H)$  can be computed in $O(|V|+\sum_{e\in E}|e|)$ time. 
 }
\end{theorem}
{

{\bf Proof}.  Consider the following algorithm. For  $i=1,.2,...,n$, select a vertex   $x_i$ of  of minimum strong degree
$s_i={s}(H_i)$ in the induced subhypergraph
$H_i=H[V_i] $  whose  vertex set is  $V_i=V-\{x_1,x_2,...,x_{i-1}\}$ and whose edge set is denoted by $E_i$.
Let  $s=\max\{s_i, i=1,2,...,n\}$. We claim that ${\hat s}(H)=s$.  Note that  ${\hat s}(H)\ge s$. We will  show  that
${\hat s}(H)\le s$.
 Now let $I=(W,F)$ be an  induced  subhypergraph of $H$ whose minimum strong degree equals 
 $ {\hat s}(H)$ and let $j, 1\le j\le n,$ be the smallest integer so that $x_j\in W$. Then $s_{I}(x_j)\le s_j={s}(H_j)\le s$, since $W\subseteq E_i$ and  and consequently the claim is proved.
 To verify the claim for time  complexity, one needs to represent $H$ as a bipartite graph $H$ as the  input of  algorithm. The details are omitted.
 $\Box$

\section{Applications in  domination theory}

For a graph $G$ on vertex set $V$ and $x\in V$ let  $N(x)$ denote the  {\it open neighborhood} of $x$, that is  the set of all  vertices adjacent to $x$, not including $x$. 
 The {\it closed neighborhood} of $x$ is $N[x]=N(x)\cup \{x\}$. The closed (open) neighborhood hypergraph of an $n$ vertex  graph $G$ is a hypergraph on the same vertices  as $G$ whose  edges are all $n$ closed (open) neighborhoods of $G$. 
 A subset of vertices $S$  in  $G$ is a  {\it dominating set} \cite{HHS}, if for every  vertex $x$ in $G$,  $N[x]\cap S\ne \emptyset$. 
$S$ is a {\it total} or {\it open  domination  set}  if, 
$N(x)\cap S\ne\emptyset$. 
  $S$  is a 2-packing (packing) , if  for any distinct pair $x,y\in S$, $N[x]$ and $N[y]$ do not intersect. $S$ is an open 
  2-packing(packing), if  for any distinct pair $x,y\in S$, $N(x)$ and $N(y])$ do not intersect
  Let $\gamma(G), \gamma^o(G),  \alpha_2(G)$ and $\alpha^0_2(G)$ denote the sizes of a smallest dominating, a smallest open domination, a  largest packing and a largest open packing, respectively,  in $G$. 
Computing $\gamma(G), \gamma^oi(G), \alpha_2(G)$ and $\alpha^o_2(G)$ are known to be NP-hard.  $\gamma(G)$ can be approximated within a factor of  $O(log(n))$ times form  its optimal solution in $O(n+m)$ time, where $n$ and  $m$ are the number of vertices and  edges of $G$. The approximation algorithm  is arising from the  approximation algorithm for the set cover problem\cite{lo}  \cite{ch}. It is  known that one can not  improve the approximation factor of $O(log(n))$ asymptotically.

Let $G$ be a graph on vertex set $V$. The {\it closed neighborhood  hypergraph},  of  $G$  is a  hypergraph on vertex set $V$ and edge set $\{N[x], x\in V\}$. The  {\it open neighborhood } hypergraph of graph $G$  is a hypergraph on the vertex set $V$ and the edge set $\{N(x), x\in V\}$.
The following summarizes basic  properties of neighborhood  hypergraphs  as they  relate  to our work.

\begin{ob}\label{o1}
{\sl Let  $H$ the closed neighborhood  hypergraph of a graph $G$ with  the vertex set $V$. 
\begin{enumerate}[(i)]

\item Let  $S\subseteq V$, then $S$ is a dominating set in $G$ if and only if $S$  is an edge cover in $H$.  

\item Let $S\subseteq V$, then $S$ is a  packing in $G$ if and if $S$ is an independent set in $H$. 

\item Let $x\in V$, then $s_H(x)\le deg(x)+1$, where $deg(x)$ is degree of $x$ in $G$.  Consequently,  ${\hat  s}(H)\le \Delta(G)+1$, where $\Delta(G)$ is the maximum degree of $G$. 

\item If $G$ is a tree and $x\in V$ is a leaf, then $s_H(x)=1$. 

\end{enumerate}
}
\end{ob}

\begin{rem}\label{r1}
{\sl Observation  \ref{o1} is valid   if $H$ is the open neighborhood  hypergraph of $G$, with the  exception that in item $(iii)$, one has $s_H(x)\le deg(x)$ and consequently  ${\hat  s}(H)\le \Delta(G)$.}
\end{rem}

By the above observation,  if we  apply the greedy algorithm in Theorem \ref{t0} to the  neighborhood  hypergraph of a graph $G$,  we obtain a dominating (total domination)  set $C$ and a packing (open packing)  $X$ so that $|C|\le  {\hat{s}}(H)|X|$.  To determine how small is $C$, we need to estimate $ {\hat s}(H)$, for the  hypergraph $H$. 
As stated above, we only know ${\hat s(H)}\le \Delta(G)+1$, where $\Delta(G)$ is the maximum degree of $G$. For trees one can get a significantly better result.

Let $T$ be a tree and let  $T_1$ be a tree which is obtained after removing all leaves of $T$. 
Then  each leaf in  $T_1$ is a support vertex in $T$ (attached to a leaf) and is called a {\it canonical support vertex} in  $T$.

Next we derive two  classical results in domination theory that were   proved first  proved in \cite{MM} and \cite{Rall}, respectively.

\begin{theorem}\label{t3}
{\sl Let $T$ be a tree on the vertex set $V$ whose  closed and open neighborhood   hypergraphs are $H$ and $H^o$, respectively.  Then, the following hold.

\begin{enumerate}[(i)]

\item ${\hat {s}}(H)={\hat m}(H)=1$ and consequently 
$\gamma(T)=\alpha_2(T)$. 

\item  ${\hat {s}}(H^o)={\hat m}(H^o)=1$ and  consequently 
$\gamma^0(T)=\alpha^o_2(T)$.

\end{enumerate} 
Moreover, the domination and packing sets can be obtained in $O(V)$ time}
\end{theorem} 
{\bf Proof.} 
We first verify that  at each iteration of the greedy algorithm in Theorem \ref{t1}  a vertex of strong degree one is detected. 
This shows ${\hat si}(H)=1$. We then apply the greedy algorithm in Theorem \ref {t0} to obtain the equality of packing and 
domination numbers. 

To prove the first the claim, note that  algorithm in Theorem \ref{t1}  can break the ties arbitrary. So assume that the algorithm selects the leaves in $T$ which   as stated in \ref{o1} have strong degree one  in $H$. Now Consider the execution of algorithm on Tree $T_1$ which is obtained after removing all leaves of $T$. If $T_1$ is empty we are done, since all vertices have already had degree one. So assume $T_1$ is not empty. 

{\bf Claim.} Let $x$ be a leaf in $T_1$, then $s_I(x)=1$, where $I$ is the closed induced neighborhood hypergraph which is obtained after removal of all leaves of $T$.

{\bf Proof of claim}. Since  $x$ is leaf in $T_1$, there is  exactly one vertex $z$ adjacent to  $x$ in $T_1$. Now Let $Y\subset V$ be the set of leaves of $T$ adjacent to $x$ (in $T$) and $N_I[y]$ denote the closed  neighborhood of $y\in Y$ in $I$ after removal of $y$.  Then,  we have $N_I[y]=x$. Additionally, note that $N_I[x]=\{x,z\}\subseteq  N_I[z]$, since $x\in N_I[z]$,  and consequently $s_I(x)=1$. \\

Coming back to the proof, now let algorithm select leaves of $T_1$, then,  delete all these leaves and continue the process with the tree obtained after  removal of these leaves. This proves ${\hat {s}}(H)=1$, consequently ${\hat {m}}(H)=1$. Now run the algorithm in 
Theorem \ref{t0} on $T$ to prove $\gamma(T)=\alpha_2(T)$.

Proof of second the claim is similar to the first and is omitted. The claim on the time complexity follows from  running times stated in Theorems \ref{t0}, \ref{t1}. $\Box$ 

\section{The gap between ${\hat m}(H)$ and ${\hat s}(H)$} 
 In the proof of Theorem \ref{t3}, we were able to effectively use ${\hat s}(H)$ instead of   ${\hat m }(H)$. However, in general this may not be possible since  
 ${\hat s}(H)$ can be much larger than ${\hat m }(H)$ as demonstrated in the following. 
 
\begin{theorem}\label{t10}
 {\sl For any integer  $n\ge 3$ there is an $n$ vertex hypergraph such  that ${\hat m}(H)=2$ and ${\hat s}(H)=n-2$. 
 }
 \end{theorem} 
 
 {\bf Proof.} Let $G$ be a graph on vertex set $V=\{v_1, v_2,...,v_n\}$   composed of a clique on vertex set 
 $\{v_2,v_3,...,v_n\}$ so that vertex $v_1$ (which is not in the clique) is adjacent to vertex $v_2$ (which is in the clique).   Now define a hypergraph $H=(V,E)$ with $E=N[v1]\cup_{i=2}^n N(v_i)$.  
 Note that     

 \begin{equation}\label{e10}
 s_H(v_1)=2 
  \end{equation} 
  since  $N[v_1]$  {and}  $N(v_2)$  are maximal edges of $H$ { containing} 
$v_1$. It is also easy to verify that  
 \begin{equation}\label{e11} 
 s_H(v_2)=n-1 \mbox { and that }  s_H(v_i)=n-2 \mbox{ for  }  i=3, 4,...,n  
 \end{equation} 
  Next note that the only  strong subset of $V$ in $H$ is  $V$ itself and thus equations \ref{e10}  and \ref{e11} imply 
 ${\hat m}(H)=s_H(v_1)=2$ as claimed for mighty degeneracy.

  Now consider the induced hypergraph $I$ on vertex set $W=V-\{v_1\}$, whose edges are obtained by removing $v_1$ from those  edges of $H$   that  contains $v_1$ (these edges  are $N[v_1]$ and $N(v_2)$). One can verify that 
  \begin{equation}\label{e12}
  s_I(v_i)=n-2 \mbox{   for   } i=2,3,...,n
\end{equation} 
  which implies ${\hat s}(H)=n-2$ as claimed. $\Box$

\section{Future Work}
This paper contains our preliminary results and we suggest several directions for future research. \\
It is not known to us yet, if  ${\hat m}(H)$  can be computed in polynomial time or not. We suspect that a variation of the algorithm in Theorem \ref{t0} can actually compute  ${\hat m}(H)$, but have not been able to prove it.\\
The connections  between the  $vc(H)$ and ${\hat m}(H)$ (${\hat s}(H)$)  needs to be explored further. Is it true that one can always be bounded by a function of the other?\\

The most recent results for approximation of $\gamma(G)$ (domination number of a graph $G$) in sparse graphs require  solving the linear programming relaxations (fractional versions) of the problem  and then rounding the solutions \cite{BU},\cite{DV}. For a recent survey see   \cite{lps}.  We suspect that proper modification of our method  in Section Four would give similar results without the need to actually solve the linear programming problems.

\end{document}